\documentclass[12pt]{article}
 \usepackage{dsfont, amssymb,amsmath,amscd,latexsym, amsthm, amsxtra,amsfonts}
 \usepackage[all]{xy}
 \newtheorem{Def}{Definition}[section]
 \newtheorem{Them}{Theorem}[section]
 \newtheorem{Prop}{Proposition}[section]

\newtheorem{theorem}{Theorem}[section]

\newtheorem{proposition}[theorem]{Proposition}

\begin{document}
 \baselineskip 18pt
 \date{}
 \title{{ { \bf Anticipating Reflected Stochastic Differential Equations }} }
\author{  Zongxia Liang \\
  {\small Department of Mathematical Sciences, Tsinghua
 University,Beijing 100084, }\\{\small    China. Email:
 zliang@math.tsinghua.edu.cn} \\
  Tusheng  Zhang\\{\small Department of Mathematics,
   University of Manchester,Oxford Road,}\\{\small
 Manchester M13, 9PL,  England, U.K. Email: tzhang@maths.man.ac.uk} }
 \maketitle
 \vskip 8pt
\begin{center}\end{center} \vskip 8pt
 {\bf Abstract: }\  In
this paper, we establish  the existence of the solutions $ (X, L)$
of reflected stochastic differential equations with possible
anticipating initial random variables. The key is to obtain some
substitution formula for Stratonovich integrals via  a uniform
convergence of the corresponding Riemann
sums.\\

\begin{center}{\bf Contents}\end{center}
\S 1. Introduction and main result\\
\S 2. Regularity of the solution $ ( X_t(x), L^x_t ) $ of Eq.(1.1) \\
\S 3.  Continuity of functionals of local times   \\
\S 4.  Moments estimates for one-point and two-point motions\\
\S 5. Uniform convergence (w.r.t. x) of the Riemann Sums
      $S_\pi(t,x)$\\
\S 6. Proof of the main result \\
\S 7. Acknowledgments\\
\S 8. References \\[15pt]
{\bf MSC}(2000): Primary
60H07, 60H10, 60J60; Secondary 60J55, 60J50. \\[10pt]
 {\bf Keywords:} Reflected  stochastic differential equation; Stratonovich
 integrals; Local times;  Reflection principle; Skorohod equations;
 Substitution formulas.

 \setcounter{equation}{0}
\section{ \bf Introduction and main results}
Let $\sigma : \Re \rightarrow \Re $ be a continuous function and $B$
be an $\Re$-valued standard Brownian motion on a complete filtered
probability space $(\Omega, {\mathcal F},\{{\mathcal F}_t\}_{t\in
[0,1]}, {\bf P})$ satisfying the usual conditions. For $x\geq  0 $,
we consider the following   stochastic differential equation on
$\Re_+=[0, +\infty )$ with reflecting boundary condition:
\begin{eqnarray}
X_t(x)=x +\int^t_0 \sigma (X_s(x))\circ dB_s +L^x_t, \ \forall \ t
\in [0,1],
\end{eqnarray}
where $\circ $ denotes the Stratonovich integral. A
 pair $ (X_t(x), L_t^x, t\in [0,1])$ is called a solution to equation (1.1)
if \\
(i)$ X_0(x)=x$, $ X_t(x)\geq 0 $ for $t\in [0,1]$,\\
(ii) $X_t(x)$, $ L_t^x $ are continuous and adapted to $ \{{\mathcal
F}_t\}_{t\in [0,1]}  $,\\
(iii) $L_t^x  $ is non-decreasing with $ L_0^x=0 $ and
\begin{eqnarray}
\int^t_0 \chi_{\{ X_s(x)=0\}}dL_s^x=L_t^x,
\end{eqnarray}
(iv) $ (X_t(x), L_t^x)$ satisfies Eq.(1.1) almost surely for every
$t\geq 0$.

 There now
exists a considerable body of literature devoted to the study of
reflected stochastic differential equations(see
\cite{s3,s4,s5,s6,s7,s9, s8} and references therein ). It is
well-known that Eq.(1.1) has a unique solution for any given initial
value $x\geq 0$ if $\sigma $ and its derivative are Lipschitz
continuous
functions.\\
\vskip 0.3cm

Now consider the following question: Does there still  exist  a pair
$ (X_t, L_t, t\in [0,1])$ to solve Eq.(1.1) if the initial value is
an arbitrary non-negative random variable $Z$ which may depends on
the whole Brownian paths ? \vskip 0.3cm The answer is not
immediately clear because one needs to deal with anticipating
stochastic integration. The main purpose of the present paper is to
give an affirmative answer to the above
 question. More precisely, our main result is the following
\begin{Them}
Assume that the function $\sigma $, its derivatives $ \sigma' $ and
$ \sigma''  $ are Lipschitz continuous, and $Z$ is a nonnegative
random variable. Then there is a pair $ (X_t(Z), L_t^Z, t\in [0,1])$
that solves the following anticipating reflected  SDE,
\begin{eqnarray}
X_t(Z)=Z +\int^t_0 \sigma (X_s(Z))\circ dB_s +L^Z_t, \ \forall \ t
\in [0,1],
\end{eqnarray}
and satisfies \\
(i)$ X_0(Z)=Z$, $ X_t(Z)\geq 0 $ for $t\in [0,1]$,\\
(ii) $X_t(Z)$, $ L^Z_t $ are continuous,\\
 (iii) $L^Z_t  $ is
non-decreasing with $ L^Z_0=0 $ and
\begin{eqnarray}
\int^t_0 \chi_{\{ X_s(Z)=0\}}dL^Z_s=L^Z_t.
\end{eqnarray}
\end{Them}
Where the stochastic integral in (1.3) is interpreted as
anticipating Stratonovich integral. Let us now recall the
definition. For any $t\in [0,1]$, let $\pi$ denote an arbitrary
partition of the interval $[0,t]$ of the form: $\pi = \{0=t_0< t_1<
\cdots < t_n=t \} $. Let $ \|\pi \|=\sup\limits_{0\leq k \leq
n-1}\{(t_{k+1}-t_k ) \}$ denote the norm of $\pi $. For a stochastic
process $ f=\{ f_s, \ s\in [0,1] \}$, we define its Riemann sums
$S_\pi(f, t) $ by
\begin{eqnarray}
S_\pi(f, t)=\sum^{n-1}_{k=0} \frac{1}{t_{k+1}-t_k}\bigg (
\int^{t_{k+1}}_{t_k} f_s ds  \bigg )( B_{t_{k+1}}-B_{t_k}).
\end{eqnarray}
We have the following
\begin{Def}
We say that a stochastic  process $ f=\{ f_s, \ s\in [0,1] \}$ such
that $ \int^1_0 \chi_{\{s\leq t \}} |f_s | ds < +\infty $ a.s. is
Stratonovich integrable if the family  $S_\pi(f, t) $ converges in
probability as $\|\pi \|\rightarrow  0 $. In this case the limit
will be called the  Stratonovich integral of the  process $ f $ on
$[0,1]$ and will be denoted by $ \int^t_0 f_s\circ dB_s $.
\end{Def}
Let us now describe our approach. To prove Theorem 1.1, the natural
idea is to replace $x$ in (1.1) by the initial random variable $Z$
and prove that the pair $(X_t(Z), L_t^Z)$ satisfies the anticipating
SDE. To achieve this, the key is to establish the following
 substitution formula
\begin{equation}
\int^t_0 \sigma (X_s(x))\circ dB_s\big |_{x=Z}=\int^t_0 \sigma
(X_s(Z))\circ dB_s
\end{equation}
for all $t\in [0,1]$.\\
To obtain (1.6), it seems that we can not apply the existing
substitution formula in the literature (see [9],[10]) because the
regularity of the solution $X_t(x)$ of (1.1) with respect to the
initial value $x$ is not good enough to satisfy the required
hypothesis. Instead, we prove (1.6) by showing
 the  uniform convergence (w.r.t. x) of
the corresponding Riemann Sums  $S_\pi(t,x)$ (Theorem 5.1 in section
5 below ). The Garsia, Rodemich and Rumsey's Lemma and   moments
estimates for one-point and two-point motions(Theorem 4.1, 4.2 in
section 4 below ) will play an important role. For the proof of
(1.4), we need to study the continuity of the random field $\int^t_0
l(X_s(x))dL_s^x$ for any continuous function $l(y)$  on $(0, \infty
)$ with compact support.

\vskip 0.3cm This paper is organized as follows. Section 2 is to
study the  regularity of  the solution $ ( X_t(x), L^x_t ) $ of
Eq.(1.1). In Section 3 we prove the continuity of functional of
local times. Section 4  is to study moments estimates for one-point
and two-point motions. In Section 5 we prove  the uniform
convergence (w.r.t.x) of the Riemann Sums
      $S_\pi(t,x)$. The proof of Theorem 1.1 will be completed in Section 6.

\setcounter{equation}{0}
\section{\bf Regularity of the solution $ ( X_t(x), L^x_t ) $ of Eq.(1.1)}
We first recall the deterministic Skorohod  problem(see\cite{s1}).
\begin{Def}
Let $ y\in \big\{f\in C( [0,1]\longrightarrow \Re ), \ f(0)\geq 0
\big \}$. We will  say that a pair $(x, k) $ of functions on $[0,1]$
is a solution of the Skorohod  problem  associated with $y$ if \\
(i) $x_t=y_t+ k_t$, $t\in [0,1]$,\\
(ii) $x_t\geq 0 $, $t\in [0,1]$,\\
(iii) $k$ is increasing, continuous, $k(0)=0 $ and satisfies
\begin{eqnarray}
k_t=\int^t_0 \chi_{\{x_s=0\}}dk_s.
\end{eqnarray}
In this case, the function $k$ is given by
\begin{eqnarray}
k_t=-\inf_{s\leq t}\{ (y_s \wedge 0)\}.
\end{eqnarray}
\end{Def}

\begin{Prop}
Assume that the function $ \sigma $ satisfies the same conditions as
in Theorem 1.1, and  $ ( X_t(x), L^x_t ) $  is a solution of
Eq.(1.1). Then there is a constant $c$ such that
\begin{eqnarray}
{\bf E}\big \{\sup_{t\leq 1}|X_t(x)-X_t(y)|^{p}  \big \} \leq \exp\{
cp^2 +cp\}|x-y|^p
\end{eqnarray}
for any $x, y \in \Re_+$ and $p\geq 1$.
\end{Prop}
{\bf Proof.}\ By H\"{o}lder inequality, we need only to prove
Proposition 2.1 for $p\geq 4$. \  Let $ a(x):= \frac{1}{2}(\sigma
\sigma')(x) $ for any $x\in \Re $, we have
\begin{eqnarray}
&& X_t(x)=x + \int^t_0 \sigma (X_s(x)) dB_s + \int^t_0 a(X_s(x)) ds
+L^x_t,\ t\in [0,1],\nonumber\\ \\
&& X_t(y)=y + \int^t_0 \sigma (X_s(y)) dB_s + \int^t_0 a(X_s(y)) ds
+L^y_t,\ t\in [0,1].\nonumber\\
\end{eqnarray}
By the reflection principle (2.2),
\begin{eqnarray}
&& L^x_t=-\inf_{s\leq t}\big \{\big (x + \int^t_0 \sigma (X_s(x))
dB_s + \int^t_0 a(X_s(x)) ds      \big )\wedge 0      \big \},\ t\in
[0,1], \nonumber\\
\\
&& L^y_t=-\inf_{s\leq t}\big \{\big (y + \int^t_0 \sigma (X_s(y))
dB_s + \int^t_0 a(X_s(y)) ds      \big )\wedge 0      \big \},\ t\in
[0,1].\nonumber\\
\end{eqnarray}
Thus,
\begin{eqnarray}
 |X_t(x)-X_t(y)|& \leq &  |x-y| + \sup_{s\leq t}\big |\int^t_0 [ \sigma
(X_s(x))- \sigma (X_s(y)) ]dB_s\big | + | L^x_t-L^y_t|\nonumber\\
&\leq & 2 |x-y| +2 \sup_{s\leq t}\big |\int^t_0 [ \sigma (X_s(x))-
\sigma (X_s(y)) ]dB_s\big |\nonumber\\
&&\quad  +2\int^t_0 | a (X_s(x))- a (X_s(y)) |ds.
\end{eqnarray}
Set $ Y_t(x, y)=\sup_{s\leq t}\{|X_s(x)-X_s(y)|\}$,
$\psi_t(x,y)=\|Y_t(x, y)\|_p $ for any $ x, y \in \Re_+ $ and $p\geq
4$. By Burkh\"o{}lder (see \cite{s12}) and H\"{o}lder inequalities,
\begin{eqnarray}
&& \big \|\sup_{s\leq t}\big |\int^t_0 [ \sigma (X_s(x))- \sigma
(X_s(y)) ]dB_s \big | \big \|_p\nonumber\\
&& \leq c p^{\frac{1}{2}}\| \big ( \int^t_0[ \sigma (X_s(x))- \sigma
(X_s(y)) ]^2ds      \big )^{\frac{1}{2}} \|_p\nonumber\\
&& \leq c p^{\frac{1}{2}}( \int^t_0 \psi_s(x,y) ^2ds      \big
)^{\frac{1}{2}},
\end{eqnarray}
where we have used  Lipschitz continuity of $ \sigma $.\\
Similarly,
\begin{eqnarray}
 \big \|\int^t_0 | a (X_s(x))- a (X_s(y)) |ds\big \|_p
  \leq c ( \int^t_0 \psi_s(x,y) ^2ds      \big
)^{\frac{1}{2}}.
\end{eqnarray}
Using (2.8), (2.9) and (2.10),we get that
\begin{eqnarray}
\psi_t(x,y) ^2\leq 12 |x-y|^2 + (12c^2p +3c^2)\int^t_0 \psi_s(x,y)
^2ds .
\end{eqnarray}
It follows from  Gronwall's lemma  and (2.11) that
\begin{eqnarray*}
{\bf E}\big \{\sup_{t\leq 1}|X_t(x)-X_t(y)|^{p}  \big \} \leq \exp\{
12cp^2 +3c^2p+3p\}|x-y|^p.
\end{eqnarray*}
Thus we complete the proof. \quad $\Box $

In view of  (2.6) and (2.7), the following is a direct consequence
of Proposition 2.1.
\begin{Prop}
Assume that the function $ \sigma $ satisfies the same conditions as
in Theorem 1.1, and $ ( X_t(x), L^x_t ) $  is a solution of
Eq.(1.1). Then there is a constant $c$ such that
\begin{eqnarray}
{\bf E}\big \{\sup_{t\leq 1}|L^x_t-L^y_t|^{p}  \big \} \leq \exp\{
cp^2 +cp\}|x-y|^p
\end{eqnarray}
for any $x, y \in \Re_+$ and $p\geq 1$.
\end{Prop}
Similar arguments lead to  the following result.
\begin{Prop}
Assume that the function $ \sigma $ satisfies the same conditions as
in Theorem 1.1, $ ( X_t(x), L^x_t ) $  is a solution of Eq.(1.1).
Then there is a constant $c$ such that
\begin{eqnarray}
&&{\bf E}\big \{\sup_{0\leq t\leq 1}|X_t(x)|^{p}  \big \} \leq
\exp\{
cp^2 +cp\}(1+x)^p,\\
&&{\bf E}\big \{\sup_{0\leq t\leq 1}|L^x_t|^{p}  \big \} \leq \exp\{
cp^2 +cp\}(1+x)^p
\end{eqnarray}
for any $x \in \Re_+$ and $p\geq 1$.
\end{Prop}
Next, we  study the  regularity of the solution $ ( X_t(x), L^x_t )
$ of Eq.(1.1) w.r.t. $t$.
\begin{Prop}
Assume that the function $ \sigma $ satisfies the same conditions as
in Theorem 1.1, $ ( X_t(x), L^x_t ) $  is a solution of Eq.(1.1).
Then,  for any $ R>0 $ and $p\geq 1 $,  there exist constants $C(p,
R)$ such that, for $s, t \in [0, 1]$,
\begin{eqnarray}
&&\sup_{0\leq x \leq R}{\bf E}\big\{ |X_t(x)-X_s(x)|^{2p}
\big \} \leq C(p, R)|t-s|^p,\\
&&\sup_{0\leq x \leq R}{\bf E}\big\{ |L^x_t-L^x_s|^{2p} \big \} \leq
C(p, R)|t-s|^p.
\end{eqnarray}
\end{Prop}
{\bf Proof.} For $ 0\leq s\leq t \leq 1 $, let $f(t)\equiv -x
-\int^t_0 \sigma (X_s(x)) dB_s - \int^t_0 a(X_s(x)) ds$. By
reflection principle,
$$ L_t^x=\sup_{0\leq u \leq t}\{f(u)\vee 0    \}       .$$
Noting that
\begin{equation}
L^x_t-L^x_s \leq \sup_{s\leq u \leq t}\{| f(u) - f(s)| \},
\end{equation}
 by Burkh\"o{}lder (see \cite{s12}) and H\"{o}lder inequalities, we have
\begin{eqnarray}
{\bf E}\{|L^x_t-L^x_s|^{2p}\}& \leq & c(p) {\bf E} \{\sup_{s\leq u
\leq t}|\int^u_s \sigma (X_v(x)) dB_v          |^{2p} \}\nonumber\\
&&+ c(p) {\bf E} \{\sup_{s\leq u \leq t}|\int^u_s a(X_v(x)) dv
|^{2p} \}\nonumber\\
&\leq & c(p) {\bf E} \big (\int^t_s \sigma (X_v(x))^2 dv \big
)^p\nonumber\\
&&  +c(p){\bf E} \big ( |\int^t_s a(X_v(x)) dv | \big
)^{2p}\nonumber\\
&\leq & c(p)( 1+R^{2p})( |t-s|^p + |t-s|^{2p})\nonumber\\
&\leq & c(p, R )|t-s|^p,
\end{eqnarray}
and this implies (2.16).
 Since
$$ | X_t(x) - X_s(x) | \leq |\int^u_s \sigma (X_v(x)) dB_v    |
+  |\int^t_s a(X_v(x)) dv |  + | L^x_t-L^x_s     |, $$ the
inequality (2.15) follows from (2.16), Burkh\"o{}lder (see
\cite{s12}) and H\"{o}lder inequalities. Thus we complete the proof
of Proposition 2.4. \quad $\Box$

 \setcounter{equation}{0}
\section{ \bf Continuity of functionals of local times }
Let $ ( X_t(x), L^x_t ) $  be a solution of Eq.(1.1). Because of
Propositions 2.1-2.2 and Proposition 2.4, we may assume that
$X_t(x), L_t^x$ are jointly continuous in $(t,x)$. Let  $ l(y) $ be
a continuous function  on $(0, \infty )$ with compact support. Put
$F(t, x)=\int^t_0 l(X_s(x))dL^x_s$ for $x\in \Re_+$. We have the
following
\begin{Prop}
Assume that the function $ \sigma $ satisfies the same conditions as
in Theorem 1.1, and $ ( X_t(x), L^x_t ) $  is a solution of
Eq.(1.1). Then the function $F(t,x)$ is jointly continuous in $(t,
x)$.
\end{Prop}
{\bf Proof.} \ Since
\begin{eqnarray}
|F(t,x)-F(s,x)|\leq \sup_{y>0}\{|l(y)|\}|L^x_t-L^x_s|,
\end{eqnarray}
the function $F(t,x)$ is continuous in $t$ uniformly with respect to
$x$ in any compact set by Proposition 2.4 and Kolmogorov's
continuity criterion( see Theorem 1.4.1 in \cite{s20} ). Thus, it
suffices to show the continuity of $ F(t,x)$ w.r.t. $x$ for any
fixed $t$. Let $x_n,x\in \Re_+ $ with  $ x_n \longrightarrow x $ as
$ n\longrightarrow +\infty $. By Propositions 2.1-2.2, and
Kolmogorov's continuity criterion( see Theorem 1.4.1 in \cite{s20}
), we have
\begin{eqnarray}
L^{x_n}_t\longrightarrow L^x_t, \quad  X_t(x_n)\longrightarrow
X_t(x),
\end{eqnarray}
uniformly in $t$, as $n\longrightarrow +\infty $. Therefore, there
exists a constant $C\geq 1$ such that for all  $n\geq 1$
\begin{eqnarray}
 L^{x_n}_t \leq C + L^x_1.
 \end{eqnarray}
Since the function $l(x)$ is bounded  and continuous, by (3.2) and
(3.3),
\begin{eqnarray}
| \int^t_0 [ l(X_s(x_n))-l(X_s(x))]dL^{x_n}_s|\longrightarrow 0
\end{eqnarray}
as $n\longrightarrow +\infty$. Because $L^{x_n}_t $ and $ L^{x}_t$
are increasing and continuous, by (3.2), the  sequence of finite
measures $ dL^{x_n}_t $ on $[0,1]$ converges weakly to the finite
measure $ dL^{x}_t $ on $[0,1]$. Therefore, for bounded continuous
function $ l(X_s(x))$ on $[0,1]$, we have
\begin{eqnarray}
\lim_{n\rightarrow \infty }\int^t_0 l( X_s(x))dL^{x_n}_s=\int^t_0 l(
X_s(x))dL^{x}_s.
\end{eqnarray}
The proof of Proposition 3.1 follows from (3.4) and (3.5).\quad
$\Box $

\setcounter{equation}{0}
\section{{\bf  Moments estimates for
one-point and two-point motions}}

For any $R >0 $ and  $ x \in [0, R] $, let $ ( X_t(x), L^x_t ) $
 be a solution of Eq.(1.1).We define  $ S_\pi (t,x) $ and $I(t,x) $
by\begin{eqnarray*}
S_\pi (t,x)&:=&S_\pi(\sigma (X_{\cdot} (x)), t),\\
I(t, x)&:=& \int^t_0\sigma (X_s (x))\circ dB_s
        =\int^t_0\sigma (X_s (x))dB_s +\frac{1}{2}
        \int^t_0 (  \sigma     \sigma' )(X_s (x)) ds.
\end{eqnarray*}
Write
\begin{eqnarray}
&&S_\pi(t,x)\nonumber\\
&=&\sum^{n-1}_{k=0} \frac{1}{t_{k+1}-t_k}\bigg (
\int^{t_{k+1}}_{t_k} \sigma ( X_s(x) ) ds  \bigg )(
B_{t_{k+1}}-B_{t_k})\nonumber\\
&=& \sum^{n-1}_{k=0}\sigma ( X_{t_k}(x) )(
B_{t_{k+1}}-B_{t_k})\nonumber\\
&+&\sum^{n-1}_{k=0} \frac{1}{t_{k+1}-t_k}\bigg (
\int^{t_{k+1}}_{t_k} (\sigma ( X_s(x) )-\sigma ( X_{t_k}(x) )) ds
\bigg )( B_{t_{k+1}}-B_{t_k}).\nonumber\\
\end{eqnarray}
By Ito's formula, for $s\geq t_k$,
\begin{eqnarray}
\sigma ( X_s(x) )&-&\sigma ( X_{t_k}(x) )\nonumber\\
&=&\int_{t_k}^s\sigma^{\prime} ( X_u(x) )\sigma ( X_u(x) )
dB_u+ \int_{t_k}^s\sigma^{\prime} ( X_u(x) ) dL_u^x \nonumber\\
&+&\frac{1}{2}\int_{t_k}^s(\sigma^{\prime})^2( X_u(x) ) \sigma (
X_u(x) ) du+\frac{1}{2}\int_{t_k}^s \sigma^{\prime\prime}( X_u(x)
)\sigma^2 ( X_u(x)) du.\nonumber\\
\end{eqnarray}
Thus we can write $S_\pi(t,x)-I(t,x )$ as follows:
\begin{equation}
S_\pi(t,x)-I(t,x )=A_{1\pi}+A_{2\pi}+A_{3\pi}+A_{4\pi},
\end{equation}
where
\begin{eqnarray*}
A_{1\pi}(x)&:=&\sum^{n-1}_{i=0} \sigma
(X_{t_i}( x))( B_{t_{i+1}}-B_{t_i})-\int^1_0 \sigma(X_s(x))dB_s,\\
A_{2\pi}(x)&:= &\sum^{n-1}_{i=0}\frac{1}{t_{i+1}-t_i} \int_{t_i}^{t_{i+1}}ds
\big (\int_{t_i}^{s}\sigma'(X_u(x))\sigma(X_u(x))dB_u\big )(B_{t_{i+1}}-B_{t_i})\\
&&- \frac{1}{2}\int_{0}^{t}\sigma'(X_s(x))\sigma(X_s(x))ds,\\
A_{3\pi}(x)&:=& \sum^{n-1}_{i=0}\frac{1}{t_{i+1}-t_i}
\int_{t_i}^{t_{i+1}}ds \big ( \int_{t_k}^s[(\sigma^{\prime})^2(
X_u(x) )
\sigma ( X_u(x) )\nonumber\\
&&+\sigma^{\prime\prime}( X_u(x) )\sigma^2 ( X_u(x)) ]du\big ) \times (B_{t_{i+1}}-B_{t_i}),\\
A_{4\pi}(x)&:=& \sum^{n-1}_{i=0} \frac{1}{t_{i+1}-t_i}
\int_{t_i}^{t_{i+1}}ds\big (
 \int_{t_k}^s\sigma^{\prime}( X_u(x) )dL_u^x\big )\times (B_{t_{i+1}}-B_{t_i}).
\end{eqnarray*}

\begin{proposition}
Assume that the function $ \sigma $ satisfies the same conditions as
in Theorem 1.1, $ ( X_t(x), L^x_t ) $  is a solution of Eq.(1.1).
Then  for any $ p\geq 2$ and $R>0$ there exist  constants $C(p, R )
$
 such that
\begin{equation}
\sup_{ x\in [0,  R]}{\bf E} \big \{|S_\pi(t,x)-I(t,x ) |^{2p}\big
\}\leq C(p, R ) \|\pi\|^{\frac{1}{2}p}.
\end{equation}
\end{proposition}

{\bf Proof.} In the sequel, we will use $c(p)$ to denote a generic
constant which depends only on $p$ and whose value may be different
from line to line. By Burkholder-Davis-Gundy inequalities, we have
\begin{eqnarray}
\{{\bf E} \{ |A_{1\pi}(x)|^{2p}\}\}^{\frac{1}{p}}&\leq & c(p) \bigg
[
{\bf E} \bigg ( \sum^{n-1}_{i=0}\int_{t_i}^{t_{i+1}}(\sigma(X_s(x))-
\sigma(X_{t_i}(x))^2ds\bigg )^p \bigg ]^{\frac{1}{p}  }\nonumber\\
&\leq & c(p)\sum^{n-1}_{i=0}\bigg \{ {\bf E}
\{\vert\int_{t_i}^{t_{i+1}}(\sigma(X_s(x))-\sigma(X_{t_i}(x))^2ds\vert^{p}
\} \bigg
\}^{\frac{1}{p}}\nonumber\\
&\leq&c(p)\sum^{n-1}_{i=0} \{\vert\int_{t_i}^{t_{i+1}}({\bf E}[(\sigma(X_s(x))
-\sigma(X_{t_i}(x))^{2p}])^{\frac{1}{p}}ds\nonumber \\
&\leq & c(p)\|\pi\|,
\end{eqnarray}
where we have used Proposition 2.4. For $p\geq 1$, by H\"{o}lder
inequality,
\begin{eqnarray}
&&{\bf E} \{ |A_{3\pi}(x)|^{2p}\}\nonumber\\
&\leq &{\bf E} \{ \bigg (\sum^{n-1}_{i=0}(\frac{1}{t_{i+1}-t_i}
 \int_{t_i}^{t_{i+1}}ds \int_{t_i}^s[(\sigma^{\prime})^2( X_u(x) )
 \sigma ( X_u(x) )\nonumber\\
 &&+\sigma^{\prime\prime}( X_u(x) )\sigma^2 ( X_u(x)) ]du)^2\bigg )^p
 \times \big (\sum^{n-1}_{i=0}(B_{t_{i+1}}-B_{t_i})^2\big )^p\}\nonumber\\
&\leq & c(p)\|\pi\|^p\sup_{\pi}{\bf E} \{\big
(\sum^{n-1}_{i=0}(B_{t_{i+1}}
-B_{t_i})^2\big )^p\}\nonumber\\
&\leq& c(p)\|\pi\|^p .
\end{eqnarray}
\begin{eqnarray}
{\bf E} \{ |A_{4\pi}(x)|^{2p}\}
&\leq &{\bf E} \{ \big (\sum^{n-1}_{i=0}(\frac{1}{t_{i+1}-t_i} \int_{t_i}^{t_{i+1}}ds
 \int_{t_i}^s\sigma^{\prime}( X_u(x) )dL_u^x)^2\big )^p\nonumber\\
&&\times \big (\sum^{n-1}_{i=0}(B_{t_{i+1}}-B_{t_i})^2\big )^p\}\nonumber\\
&\leq & c(p){\bf E} \{ \big (\sum^{n-1}_{i=0}(L_{t_{i+1}}^x-L_{t_i}^x)^2\big )^p\times
 \big (\sum^{n-1}_{i=0}(B_{t_{i+1}}-B_{t_i})^2\big )^p\}\nonumber\\
&\leq & c(p){\bf E} \{ \big (\sup_i(L_{t_{i+1}}^x-L_{t_i}^x)L_1^x\big )^p\times
\big (\sum^{n-1}_{i=0}(B_{t_{i+1}}-B_{t_i})^2\big )^p\}\nonumber\\
&\leq & c(p)({\bf E}\{\big (\sup_i(L_{t_{i+1}}^x-L_{t_i}^x)\big
)^{3p}\})^{\frac{1}{3}}
 ({\bf E}\{\big (L_{1}^x\big )^{3p}\})^{\frac{1}{3}}\nonumber\\
&&\times {\bf E} \{\big (\sum^{n-1}_{i=0}(B_{t_{i+1}}-B_{t_i})^2\big )^p\}\nonumber\\
&\leq& c(p)\|\pi\|^{\frac{p}{2}},
\end{eqnarray}
where Proposition 2.4 and Kolmogorov's continuity criterion( see
Theorem 1.4.1 in \cite{s20} ) were used in the last inequality.
Using Fubini Theorem,  $A_{2\pi}$ can be further written as
\begin{equation}
A_{2\pi}(x)=A_{2\pi}^{(1)}(x)+A_{2\pi}^{(2)}(x)+A_{2\pi}^{(3)}(x),
\end{equation}
where
\begin{eqnarray*}
A_{2\pi}^{(1)}(x)&:=& \sum^{n-1}_{i=0}\frac{1}{t_{i+1}-t_i}
\int_{t_i}^{t_{i+1}}(t_{i+1}-u) \bigg (\sigma^{\prime}( X_u(x))
\sigma ( X_u(x) )\\
&&-\sigma^{\prime}( X_{t_i}(x))\sigma ( X_{t_i}(x) \bigg )du,\\
A_{2\pi}^{(2)}(x)&:=&
-\frac{1}{2}\sum^{n-1}_{i=0}\int_{t_i}^{t_{i+1}} \big
(\sigma^{\prime}( X_u(x))\sigma ( X_u(x) )-\sigma^{\prime}(
X_{t_i}(x))
\sigma ( X_{t_i}(x) \big )du,\\
A_{2\pi}^{(3)}(x)&:= &\sum^{n-1}_{i=0}\{ \frac{1}{t_{i+1}-t_i} \big
(\int_{t_i}^{t_{i+1}}(t_{i+1}-u)\sigma'(X_u(x))\sigma(X_u(x))
dB_u\big )(B_{t_{i+1}}-B_{t_i})\\
&&-\frac{1}{t_{i+1}-t_i} \big (\int_{t_i}^{t_{i+1}}(t_{i+1}-u)
\sigma'(X_u(x))\sigma(X_u(x))du\big )\}.
\end{eqnarray*}
Since $\sigma^{\prime}\sigma$ is Lipschitz continuous, it follows
that
\begin{eqnarray}
\{{\bf E} \{ |A_{2\pi}^{(1)}(x)|^{p}\}\}^{\frac{1}{p}}&\leq &
\sum^{n-1}_{i=0}\int_{t_i}^{t_{i+1}} \big \{{\bf E}
\{\|\sigma^{\prime}( X_u(x))
\sigma ( X_u(x) )\nonumber\\
&&-\sigma^{\prime}( X_{t_i}(x))\sigma ( X_{t_i}(x) \|^p\}
\big \}^{\frac{1}{p}}du\nonumber\\
&\leq & c\sum^{n-1}_{i=0}\int_{t_i}^{t_{i+1}} \{{\bf E}
\{\| X_u(x)-X_{t_i}(x) \|^p\}\}^{\frac{1}{p}}du\nonumber\\
&\leq & c\|\pi\|^{\frac{1}{2}}.
\end{eqnarray}
Similar arguments lead to
\begin{equation}
\{{\bf E} \{ |A_{2\pi}^{(2)}(x)|^{p}\}\}^{\frac{1}{p}}\leq
c\|\pi\|^{\frac{1}{2}}.
\end{equation}
Noting that $A_{2\pi}^{(3)}$ is a  martingale. Using
Burkholder-Davis-Gundy inequalities, we obtain that
\begin{eqnarray}
&&\{{\bf E} \{ |A_{2\pi}^{(3)}(x)|^{2p}\}\}^{\frac{1}{p}}\nonumber\\
&\leq & c(p) \bigg [ {\bf E} \bigg ( \sum^{n-1}_{i=0}\{
\frac{1}{t_{i+1}-t_i} \big
(\int_{t_i}^{t_{i+1}}(t_{i+1}-u)\sigma'(X_u(x))\sigma(X_u(x))dB_u\big
)
(B_{t_{i+1}}-B_{t_i})\nonumber\\
&&-\frac{1}{t_{i+1}-t_i} \big (\int_{t_i}^{t_{i+1}}(t_{i+1}-u)
\sigma'(X_u(x))\sigma(X_u(x))du\big )\}^2\bigg )^p \bigg
]^{\frac{1}{p}}
\nonumber\\
&\leq& c(p)\sum^{n-1}_{i=0}\bigg [ {\bf E}\{ \frac{1}{t_{i+1}-t_i}
\big (\int_{t_i}^{t_{i+1}}(t_{i+1}-u)
\sigma'(X_u(x))\sigma(X_u(x))dB_u\big )(B_{t_{i+1}}-B_{t_i})\nonumber\\
&&-\frac{1}{t_{i+1}-t_i} \big (\int_{t_i}^{t_{i+1}}(t_{i+1}-u)
\sigma'(X_u(x))\sigma(X_u(x))du\big )\}^{2p} \bigg ]^{\frac{1}{p}}\nonumber\\
&\leq& c(p)\sum^{n-1}_{i=0}\bigg \{ \big ({\bf E}\{ \big
|\frac{1}{t_{i+1}-t_i} \int_{t_i}^{t_{i+1}}
(t_{i+1}-u)\sigma'(X_u(x))\sigma(X_u(x))dB_u\big )(B_{t_{i+1}}
-B_{t_i})\|^{2p}\big )^{\frac{1}{p}}\nonumber\\
&&+\big ({\bf E}\{ \big |\frac{1}{t_{i+1}-t_i} \big
(\int_{t_i}^{t_{i+1}}(t_{i+1}-u)\sigma'(X_u(x))
\sigma(X_u(x))du\big |\}^{2p}\big )^{\frac{1}{p}}\bigg \}\nonumber\\
&\leq & c(p)\|\pi\|,
\end{eqnarray}
where H\"o{}lder inequality was used for the last inequality.
Combining the estimates for $A_{2\pi}^{(1)}$, $A_{2\pi}^{(2)}$ and
$A_{2\pi}^{(3)}$  together, we deduce that
\begin{equation}
\{{\bf E} \{ |A_{2\pi}(x)|^{2p}\}\}^{\frac{1}{p}}\leq
c\|\pi\|^{\frac{1}{2}}.
\end{equation}
Now (4.4) follows from (4.5), (4.6), (4.7) and (4.12). The proof is
complete. \quad $\Box$ \vskip 0.4cm

 Next result is the moment estimates for
the two point motions.
\begin{proposition}
Assume that the function $ \sigma $ satisfies the same conditions as
in Theorem 1.1, $ ( X_t(x), L^x_t ) $  is a solution of Eq.(1.1).
Then  for any $ p\geq 2$ and $R>0$ there exists  constants $C(p, R )
$ and $\beta$, independent of the partition $\pi$,
 such that
\begin{equation}
{\bf E} \big \{\sup_{ t\in [0,  1]}|S_\pi(t,x)-S_\pi(t,y) |^{p}\big
\}\leq C(p, R ) |x-y|^{\frac{1}{2}p},
\end{equation}
for all $x, y\in [0, R]$.
\end{proposition}
{\bf Proof.} Similarly as (4.3), write
\begin{equation}
S_\pi(t,x)-S_\pi(t,y)
=A_{1\pi}(x,y)+A_{2\pi}(x,y)+A_{3\pi}(x,y)+A_{4\pi}(x,y),
\end{equation}
where
\begin{eqnarray*}
A_{1\pi}(x,y)&:=&\sum^{n-1}_{i=0}\big ( \sigma (X_{t_i}( x))- \sigma
(X_{t_i}( y))\big )( B_{t_{i+1}}-B_{t_i}),\\
A_{2\pi}(x,y)&:= &\sum^{n-1}_{i=0}\big \{\frac{1}{t_{i+1}-t_i}
\bigg ( \int_{t_i}^{t_{i+1}}ds \int_{t_i}^{s}\big (\sigma'(X_u(x))\sigma(X_u(x))\\
&&-\sigma'(X_u(y))\sigma(X_u(y))\big )dB_u\bigg )\times(B_{t_{i+1}}-B_{t_i})\big \},\\
A_{3\pi}(x,y)&:=& \sum^{n-1}_{i=0}\bigg \{\frac{1}{t_{i+1}-t_i}
 \int_{t_i}^{t_{i+1}}ds\big \{ \int_{t_k}^s[(\sigma^{\prime})^2( X_u(x) )
 \sigma ( X_u(x) )\\
 &&+\sigma^{\prime\prime}( X_u(x) )\sigma^2 ( X_u(x))-(\sigma^{\prime})^2( X_u(y) )
 \sigma ( X_u(y) )\\
 &&-
\sigma^{\prime\prime}( X_u(y) )\sigma^2 ( X_u(y)) ]du\big\}
\times (B_{t_{i+1}}-B_{t_i})\bigg \},\\
A_{4\pi}(x,y)&:=& \sum^{n-1}_{i=0}\frac{1}{t_{i+1}-t_i}
 \int_{t_i}^{t_{i+1}}ds\big\{ \int_{t_i}^s\sigma^{\prime}( X_u(x) )dL_u^x
 \\
 &&-\int_{t_i}^s\sigma^{\prime}( X_u(y) )dL_u^y\big \}(B_{t_{i+1}}-B_{t_i})\\
&=& \sigma^{\prime}(0)\sum^{n-1}_{i=0}\frac{1}{t_{i+1}-t_i}
\int_{t_i}^{t_{i+1}}ds\big\{
\int_{t_i}^sdL_u^x-\int_{t_i}^sdL_u^y\big \} (B_{t_{i+1}}-B_{t_i}).
\end{eqnarray*}
By Burkholder-Davis-Gundy inequalities, the Lipschitz continuity of
$\sigma$
 and (2.3), it follows easily that
\begin{equation}
{\bf E} \big \{\sup_{ t\in [0,  1]}|A_{1\pi}(x,y) |^{p}\big \}\leq
C(p, R ) |x-y|^{p},
\end{equation}
By virtue of (2.3) and the Lipschitz continuity of $\sigma'\sigma$,
\begin{eqnarray}
&&\bigg ({\bf E} \big \{\sup_{ t\in [0,  1]}|A_{2\pi}(x,y)
|^{p}\big \}\bigg )^{\frac{1}{p}}\nonumber\\
&\leq &\sum^{n-1}_{i=0}\big \{\frac{1}{t_{i+1}-t_i}
\int_{t_i}^{t_{i+1}}ds\bigg ( {\bf E}
\big \{\vert \int_{t_i}^{s}\big (\sigma'(X_u(x))\sigma(X_u(x))\nonumber\\
&&-\sigma'(X_u(y))\sigma(X_u(y))\big )dB_u\bigg )
\times(B_{t_{i+1}}-B_{t_i})\vert^p\big \}\bigg )^{\frac{1}{p}}\big\}\nonumber\\
&\leq &\sum^{n-1}_{i=0}\big \{\frac{1}{t_{i+1}-t_i}
\int_{t_i}^{t_{i+1}}ds\bigg ( {\bf E} \big \{\vert
\int_{t_i}^{s}\big (\sigma'(X_u(x))\sigma(X_u(x))\nonumber\\
&&-\sigma'(X_u(y))\sigma(X_u(y))\big )dB_u\vert^{2p}\bigg
)^{\frac{1}{2p}}\bigg ( {\bf E} \big \{\vert (B_{t_{i+1}}
-B_{t_i})\vert^{2p}\big \}\bigg )^{\frac{1}{2p}}\big\}\nonumber\\
&\leq &\sum^{n-1}_{i=0}\big \{\frac{1}{t_{i+1}-t_i}
\int_{t_i}^{t_{i+1}}ds\bigg (\int_{t_i}^{s} \big \{{\bf E}[\vert
X_u(x))-X_u(y)\vert^{2p}]
\big\}^{\frac{1}{p}}du\bigg )^{\frac{1}{2}}\nonumber\\
&&\times (t_{i+1}-t_i)^{\frac{1}{2}} \big\}\nonumber\\
 &\leq & c(p)|x-y|.
\end{eqnarray}
By a similar, but simpler argument, we also found that
\begin{equation}
\bigg ({\bf E} \big \{\sup_{ t\in [0,  1]}|A_{3\pi}(x,y) |^{p}\big
\}\bigg )^{\frac{1}{p}}\leq c(p)|x-y|.
\end{equation}
Observe that
\begin{eqnarray}
(A_{4\pi}(x,y))^2&\leq &c\sum^{n-1}_{i=0}\frac{1}{t_{i+1}-t_i}
\int_{t_i}^{t_{i+1}}ds\big ( L_s^x-L_{t_i}^x-L_s^y+L_{t_i}^x\big )^2\nonumber\\
&&\times \sum^{n-1}_{i=0}(B_{t_{i+1}}-B_{t_i})^2\nonumber\\
&\leq &c(\sup_{0\leq s\leq 1}\vert L_s^x-L_s^y\vert)\sum^{n-1}_{i=0}
\frac{1}{t_{i+1}-t_i} \int_{t_i}^{t_{i+1}}\big (\vert  L_s^x-L_{t_i}^x\vert
 +\vert L_s^y-L_{t_i}^x\vert \big )ds\nonumber\\
&&\times \sum^{n-1}_{i=0}(B_{t_{i+1}}-B_{t_i})^2\nonumber\\
&\leq &c(\sup_{0\leq s\leq 1}\vert L_s^x-L_s^y\vert)\sum^{n-1}_{i=0}
\big (\vert  L_{t_{i+1}}^x-L_{t_i}^x\vert +\vert L_{t_{i+1}}^y-L_{t_i}^x\vert \big )\nonumber\\
&&\times \sum^{n-1}_{i=0}(B_{t_{i+1}}-B_{t_i})^2.
\end{eqnarray}
It follows that
\begin{eqnarray}
&&{\bf E} \big \{\sup_{ t\in [0,  1]}|A_{4\pi}(x,y)
|^{2p}\big \}\nonumber\\
&\leq &c\big ({\bf E} \big \{(\sup_{0\leq s\leq 1} \vert
L_s^x-L_s^y\vert^{3p}\big \}\big )^{\frac{1}{3}}
\big ({\bf E} \big \{(  L_{1}^x+L_{1}^y)^{3p}\big \})^{\frac{1}{3}}\nonumber\\
&&\times \big ({\bf E} \big \{\big (\sum^{n-1}_{i=0}
(B_{t_{i+1}}-B_{t_i})^2)\big )^{3p}\big
\}\big )^{\frac{1}{3}}\nonumber\\
&\leq & c(p, R) |x-y|^{\frac{1}{2}p}.
\end{eqnarray}
Putting together above estimates (4.15)-(4.17) and (4.19), we arrive
at (4.13). \quad $\Box$

 The following result can be proved similarly
as Proposition 4.2.
\begin{proposition}
Assume that the function $ \sigma $ satisfies the same conditions as
in Theorem 1.1, $ ( X_t(x), L^x_t ) $  is a solution of Eq.(1.1).
Then  for any $ p\geq 2$ and $R>0$ there exists  a constant $C(p, R
)$,
 such that
\begin{equation}
{\bf E} \big \{\sup_{ t\in [0,  1]}|I(t,x)-I(t,y) |^{p}\big \}\leq
C(p, R ) |x-y|^{p}
\end{equation}
for all $x, y\in [0, R]$.
\end{proposition}
\setcounter{equation}{0}
\section{{\bf Uniform convergence of the Riemann sums}}
\begin{theorem}
Assume that the function $ \sigma $ satisfies the same conditions as
in Theorem 1.1, $ ( X_t(x), L^x_t ) $  is a solution of Eq.(1.1).
Then  for any $ p\geq 2$ and $R>0$,
\begin{equation}
\lim_{\|\pi\|\rightarrow 0}{\bf E} \big \{\sup_{ x\in [0,
R]}|S_\pi(t,x)-I(t,x ) |^{2p}\big \}=0.
\end{equation}
\end{theorem}
{\bf Proof}. By Propositions 4.2 and 4.3, it follows from
Garsia-Rodemich and Rumsey's Lemma (cf.\cite{s14})that there exists
a constant $\beta_0>0$, independent of $\pi$, such that for $x, y\in
[0,R]$,
\begin{equation}
|S_\pi(t,x)-S_\pi(t,y)|\leq K_{\pi}(\omega)|x-y|^{\beta_0},
\end{equation}
\begin{equation}
|I(t,x)-I(t,y)|\leq K(\omega)|x-y|^{\beta_0},
\end{equation}
where $ K_{\pi}(\omega), K(\omega)$ are random variables that
satisfy
\begin{equation}
\sup_{\pi}{\bf E}\{|K_{\pi}|^p\}<\infty,\quad {\bf
E}\{|K|^p\}<\infty,
\end{equation}
for all $p\geq 1$. This is possible because the constant in
Propositions 4.2 is independent of $\pi$. Thus, given any
$\varepsilon>0$, there is $\delta>0$ such that
\begin{equation}
{\bf E}\{\sup_{x\in [0,R]}\sup_{y: |y-x|\leq \delta}
|S_\pi(t,y)-S_\pi(t,x)|^p\}\leq \varepsilon,
\end{equation}
\begin{equation}
{\bf E}\{\sup_{x\in [0,R]}\sup_{y: |y-x|\leq
\delta}|I(t,y)-I(t,x)|^p\}\leq \varepsilon.
\end{equation}
On the other hand, for any  $\delta>0$, we can find $x_1, ..., x_m$
such that
$$ [0,R]\subset \cup_{i=1}^m B(x_i, \delta).$$
Consequently,
\begin{eqnarray}
\sup_{ x\in [0,  R]}|S_\pi(t,x)&-&I(t,x )
|^{p}\nonumber\\
&\leq & \sup_{x_i}\sup_{y\in B(x_i,\delta)}|S_\pi(t,y)-I(t,y )|^{p}\nonumber\\
&\leq& c(p)\sup_{x_i}\sup_{y\in
B(x_i,\delta)}|S_\pi(t,y)-S_\pi(t,x_i)|^{p}
\nonumber\\
&&+c(p)\sum_{i=1}^m|S_\pi(t,x_i)-I(t,x_i)|^{p}\nonumber\\
&&+c(p)\sup_{x_i}\sup_{y\in B(x_i,\delta)}|I(t,y)-I(t,x_i)|^{2p}.
\end{eqnarray}
By virtue of (5.5) and (5.6), this implies that
\begin{eqnarray}
{\bf E}\{\sup_{ x\in [0,  R]}|S_\pi(t,x)&-&I(t,x )
|^{p}\}\nonumber\\
&\leq& 2c(p)\varepsilon+c(p){\bf
E}\{\sum_{i=1}^m|S_\pi(t,x_i)-I(t,x_i)|^p\}.
\end{eqnarray}
Let first $\|\pi\|\rightarrow 0$ and then $\varepsilon\rightarrow 0$
to get (5.1) by Proposition 4.1 and (5.8). Thus we complete the
proof. \quad $\Box$

\setcounter{equation}{0}
\section{ \bf Proof of Theorem 1.1 }
We will prove that $(X_t:=X_t(Z), L_t^Z)$ solves the anticipating
reflected SDE (1.3). Let $l(y)$ be any given continuous function on
$(0,\infty)$ with compact support. Since $ ( X_t(x), L^x_t ) $ is a
solution of Eq.(1.1),  $ F_t(x)=\int_0^tl(X_s(x))dL_s^x=0 $ for all
$\omega$ in  some measurable set  $\Omega_{t,x}$ with ${\bf P}
(\Omega_{t,x} )=1$. By continuity of $ F_t(x) $ proved in
Proposition 3.1,  we have
\begin{eqnarray*}
{\bf P}\big ( \omega\ : F_t(x, \omega)=0 \quad \mbox{for all
$(t,x)\in [0,1]\times [0, +\infty )$} \big )=1.
\end{eqnarray*}
Hence
\begin{eqnarray}
{\bf P}\big ( \omega\ : F_t(Z, \omega)=0 \quad \mbox{for all $t\in
[0,1]$} \big )=1,
\end{eqnarray}
for any non-negative random variable $Z$.  This implies that
$$L_t^Z=\int_0^t\chi_{\{X_s(Z)=0\}}dL_s^Z$$
Next we prove (1.6). By Theorem 5.1, the following  holds almost
surely on $\{\omega; Z(\omega)\leq M\}$
$$ \int^t_0 \sigma (X_s(x))\circ dB_s\big |_{x=Z}\chi_{\{\omega; Z(\omega)\leq M\}}$$
$$=\lim_{\|\pi\|\rightarrow 0}S_\pi(t,x)\big |_{x=Z}\chi_{\{\omega; Z(\omega)\leq M\}}$$
$$= \int^t_0 \sigma (X_s(Z))\circ dB_s\chi_{\{\omega; Z(\omega)\leq M\}}.$$
Letting $M\rightarrow \infty$ we obtain the substitution formula
(1.6), and therefore prove the Theorem. \quad \quad $\Box$
\vskip 6pt

{\bf Acknowledgements.} This work is supported by NSFC and SRF for
ROCS, SEM. The author would like to thank both for their generous
financial support.

\end{document}